\documentclass[11pt]{amsart}
\usepackage{txfonts}
\usepackage{mathrsfs}
\usepackage{amsmath}
\usepackage[pagewise]{lineno}
\allowdisplaybreaks[4]
\usepackage{}
\usepackage{xypic}
\usepackage{amsfonts}
\usepackage{amssymb}
\usepackage{bbm}

\usepackage[colorlinks,
            linkcolor=blue,
           anchorcolor=blue,
            citecolor=blue
            ]{hyperref}



\newtheorem{thm}{Theorem}[section]
\newtheorem{cor}[thm]{Corollary}
\newtheorem{lem}[thm]{Lemma}

\newtheorem{prop}[thm]{Proposition}
\newtheorem{rem}[thm]{Remark}
\numberwithin{equation}{section}

\newcommand{\be}{\begin{equation}}
\newcommand{\ee}{\end{equation}}
\newcommand{\bes}{\begin{eqnarray}}
\newcommand{\ees}{\end{eqnarray}}
\newcommand{\bess}{\begin{eqnarray*}}
\newcommand{\eess}{\end{eqnarray*}}

\newcommand{\bali}{\begin{align}}
\newcommand{\eali}{\end{align}}

\begin{document}
\title[The Grothendieck algebras of smash product semisimple Hopf algebras]{The Grothendieck algebras of certain smash product semisimple Hopf algebras}
%\thanks{This work was supported by }
\author{Zhihua Wang}
\address{Z. Wang\newline Department of Mathematics, Taizhou University,
Taizhou 225300, China}
\email{mailzhihua@126.com}
\author{Gongxiang Liu}
\address{G. Liu\newline Department of Mathematics, Nanjing University, Nanjing 210093, China}
\email{gxliu@nju.edu.cn}
\author{Libin Li}
\address{L. Li\newline School of Mathematical Science, Yangzhou University, Yangzhou 225002, China}
\email{lbli@yzu.edu.cn}
\date{}

\begin{abstract}
Let $H$ be a semisimple Hopf algebra over an algebraically closed field $\mathbbm{k}$ of characteristic $p>\dim_{\mathbbm{k}}(H)^{1/2}$  and $p\nmid 2\dim_{\mathbbm{k}}(H)$. In this paper, we consider the smash product semisimple Hopf algebra $H\#\mathbbm{k}G$, where $G$ is a cyclic group of order $n:=2\dim_{\mathbbm{k}}(H)$. Using irreducible representations of $H$ and those of $\mathbbm{k}G$, we determine all non-isomorphic irreducible representations of $H\#\mathbbm{k}G$. There is a close relationship between the Grothendieck algebra $(G_0(H\#\mathbbm{k}G)\otimes_{\mathbb{Z}}\mathbbm{k},*)$ of $H\#\mathbbm{k}G$ and the Grothendieck algebra $(G_0(H)\otimes_{\mathbb{Z}}\mathbbm{k},*)$ of $H$. To establish this connection, we endow with a new multiplication operator $\star$ on $G_0(H)\otimes_{\mathbb{Z}}\mathbbm{k}$ and show that the Grothendieck algebra $(G_0(H\#\mathbbm{k}G)\otimes_{\mathbb{Z}}\mathbbm{k},\ast)$ is isomorphic to the direct sum of $(G_0(H)\otimes_{\mathbb{Z}}\mathbbm{k},*)^{\oplus\frac{n}{2}}$ and $(G_0(H)\otimes_{\mathbb{Z}}\mathbbm{k},\star)^{\oplus\frac{n}{2}}$.
\end{abstract}
\maketitle
\section{\bf Introduction}
The Grothendieck rings of finite dimensional semisimple or cosemisimple Hopf algebras have been studied by Nichols and Richmond \cite{NR}, Nikshych \cite{Ni}, Kashina \cite{Ka}, Yang \cite{CYW, YY}, etc. For a finite dimensional semisimple Hopf algebra $H$, the category Rep$(H)$ of finite dimensional representations of $H$ is a fusion category. As an important invariant of Rep$(H)$, the Grothendieck ring $G_0(H)$ of $H$ reveals
the decompositions of tensor product of irreducible representations into irreducibles. Hence the Grothendieck ring $G_0(H)$ can be used to study
the fusion category Rep$(H)$. For instance, the knowledge of the structure of the Grothendieck ring $G_0(H)$ allows to determine all fusion subcategories of Rep$(H)$, which correspond to the so-called based subrings of $G_0(H)$.

For a semisimple Hopf algebra $H$ with antipode $S$ over a field $\mathbbm{k}$, it is known that $S^2$ is an inner automorphism of $H$ (see \cite{Lo}). Here an inner automorphism
is understood to be the conjugation by an invertible element of $H$. If the ground field $\mathbbm{k}$ has positive characteristic $p$, whether or not $S^2$ can be given by conjugation with a group-like element is not completely solved (this problem is closely related to the Kaplansky's fifth conjecture).
However, such a Hopf algebra $H$ can be embedded into another finite dimensional Hopf algebra $H\#\mathbbm{k}G$, namely, the smash product of $H$ and a group algebra $\mathbbm{k}G$, in which the square of the antipode is the conjugation with a group-like element. We refer to \cite{LL, Mo, So} for such Hopf algebras and related researches.

If $H$ is a semisimple involutory Hopf algebra, namely, a semisimple Hopf algebra with $S^2=id$, the smash product Hopf algebra $H\#\mathbbm{k}G$ considered here is nothing but the usual tensor product Hopf algebra $H\otimes \mathbbm{k}G$.
In this case, the representations of $H\otimes \mathbbm{k}G$ can be stemmed directly from the representations of $H$ and those of $\mathbbm{k}G$. Also, the Grothendieck algebra of $H\otimes \mathbbm{k}G$ is the usual tensor product of the Grothendieck algebra of $H$ and that of $\mathbbm{k}G$. However, if $H$ is not necessarily involutory (although the Kaplansky's fifth conjecture states that a semisimple Hopf algebra is necessarily  involutory), the relationship between the Grothendieck algebra of $H\# \mathbbm{k}G$ and that of $H$ is not clear.

The purpose of this paper is to study representations of the smash product semisimple Hopf algebra $H\# \mathbbm{k}G$ and to establish a relationship between the Grothendieck algebra of $H\# \mathbbm{k}G$ and that of $H$, where $H$ is a semisimple Hopf algebra over a field $\mathbbm{k}$ of characteristic $p>\dim_{\mathbbm{k}}(H)^{1/2}$, $G$ is a cyclic group of order $2\dim_{\mathbbm{k}}(H)$ and $p\nmid2\dim_{\mathbbm{k}}(H)$. It is worthy mentioning that such a Hopf algebra $H$ is not necessarily involutory unless the characteristic $p$ is larger than a certain number (see \cite{So,EG}).

The paper is organized as follows: In Section 2, we present some properties of a special element $\mathbf{v}$ of the semisimple Hopf algebra $H$. Such an element is used later to describe representations of the smash product semisimple Hopf algebra $H\# \mathbbm{k}G$. In Section 3, using irreducible representations of $H$ and those of $\mathbbm{k}G$ we are able to determine all non-isomorphic irreducible representations of $H\#\mathbbm{k}G$. We also describe the dual of these irreducible representations of
$H\# \mathbbm{k}G$ in this section.
In Section 4, we endow with a new multiplication operator $\star$ on the Grothendieck algebra $G_0(H)\otimes_{\mathbb{Z}}\mathbbm{k}$ so as to obtain a new algebra $(G_0(H)\otimes_{\mathbb{Z}}\mathbbm{k},\star)$. This algebra is nothing but the usual Grothendieck algebra $(G_0(H)\otimes_{\mathbb{Z}}\mathbbm{k},\ast)$ if $H$ is involutory.
We show that the Grothendieck algebra $(G_0(H\#\mathbbm{k}G)\otimes_{\mathbb{Z}}\mathbbm{k},\ast)$
%has the following direct sum decomposition as subalgebras:
%$$(G_0(H\#\mathbbm{k}G)\otimes_{\mathbb{Z}}\mathbbm{k},\ast)\cong(G_0(H)\otimes_{\mathbb{Z}}\mathbbm{k},*)^{\oplus\frac{n}{2}}\bigoplus(G_0(H)\otimes_{\mathbb{Z}}\mathbbm{k},\star)^{\oplus\frac{n}{2}}.$$
is isomorphic to a direct sum of $(G_0(H)\otimes_{\mathbb{Z}}\mathbbm{k},*)^{\oplus\frac{n}{2}}$ and $(G_0(H)\otimes_{\mathbb{Z}}\mathbbm{k},\star)^{\oplus\frac{n}{2}}$. This reveals a relationship between the Grothendieck algebra of $H\# \mathbbm{k}G$ and that of $H$.
Moreover, we find a fusion subcategory $\mathcal{C}$ of Rep$(H\# \mathbbm{k}G)$ whose Grothendieck algebra $(G_0(\mathcal{C})\otimes_{\mathbb{Z}}\mathbbm{k},\ast)$ happens to be the direct sum of $(G_0(H)\otimes_{\mathbb{Z}}\mathbbm{k},\ast)$ and $(G_0(H)\otimes_{\mathbb{Z}}\mathbbm{k},\star)$. In view of this, the Grothendieck algebra $(G_0(H\#\mathbbm{k}G)\otimes_{\mathbb{Z}}\mathbbm{k},\ast)$ is isomorphic to $(G_0(\mathcal{C})\otimes_{\mathbb{Z}}\mathbbm{k},\ast)^{\oplus\frac{n}{2}}$.

\section{\bf Preliminaries}
Throughout, $H$ is a finite dimensional semisimple Hopf algebra over an algebraically closed field $\mathbbm{k}$ of positive characteristic $p>\dim_{\mathbbm{k}}(H)^{1/2}$ and $p\nmid2\dim_{\mathbbm{k}}(H)$.
We denote $\{V_i\mid0\leq i\leq m-1\}$ the set of all non-isomorphic finite dimensional simple $H$-modules, where $V_0$ is fixed to be the trivial $H$-module $\mathbbm{k}$. The condition $p>\dim_{\mathbbm{k}}(H)^{1/2}$ is used to make sure that  $\dim_{\mathbbm{k}}(V_i)\neq0$ in $\mathbbm{k}$ for $0\leq i\leq m-1$. Indeed,
$p^2>\dim_{\mathbbm{k}}(H)=\sum_{i=0}^{m-1}\dim_{\mathbbm{k}}(V_i)^2\geq\dim_{\mathbbm{k}}(V_i)^2$ implies that $p>\dim_{\mathbbm{k}}(V_i)$, which turns out that $\dim_{\mathbbm{k}}(V_i)\neq0$ in $\mathbbm{k}$. The character afforded by simple $H$-module $V_i$ is denoted by $\chi_i$ for $0\leq i\leq m-1$ and the character afforded by left regular module $H$ is denoted by $\chi_H$.
We denote $\{e_i\mid0\leq i\leq m-1\}$ the set of all central primitive idempotents of $H$, where the central idempotent $e_i$ acts as the identity on $V_i$ and annihilates $V_j$ for $j\neq i$.

As a Hopf algebra, $H$ has a counit $\varepsilon$, antipode $S$ and comultiplication $\Delta$, where the comultiplication $\Delta(h)$ will be written as $\Delta(h)=h_{(1)}\otimes h_{(2)}$ for $h\in H$, here we omit the summation sign. We choose a left integral $\Lambda$ in $H$ and a right integral $\lambda$ in $H^*$ such that $\lambda(\Lambda)=1$.
We denote $\mathbf{u}:=S(\Lambda_{(2)})\Lambda_{(1)}$.
We refer to \cite{Mon} for basic theory of Hopf algebras.

For the Hopf algebra $H$, there is a formula for $S^2$ as follows (see \cite{WLL1}):  $$S^2(h)=\mathbf{u}h\mathbf{u}^{-1}\ \text{for}\ h\in H.$$
For the element $\mathbf{u}$, we have the following result (see \cite[Proposition 3.3]{WLL1}):

\begin{prop}\label{prop3.20}
The element $\mathbf{u}=S(\Lambda_{(2)})\Lambda_{(1)}$ satisfies the the following properties:
\begin{enumerate}
  \item $\mathbf{u}=\chi_{H}(\Lambda_{(1)})S(\Lambda_{(2)})$.
  \item $\Lambda_{(1)}\mathbf{u}^{-1} S(\Lambda_{(2)})=1.$
  \item $\lambda(e_i)=\dim_{\mathbbm{k}}(V_i)\chi_i(\mathbf{u}^{-1})$.
  \item $\mathbf{u}S(\mathbf{u})=S(\mathbf{u})\mathbf{u}=\varepsilon(\Lambda)\sum_{i=0}^{m-1}\frac{\dim_{\mathbbm{k}}(V_i)^2}{\lambda(e_i)}e_i$.
  \item $S(\mathbf{u}^{-1})\mathbf{u}=\mathbf{u}S(\mathbf{u}^{-1})$, which is the distinguished group-like element of $H$.
\end{enumerate}
\end{prop}

We denote \begin{equation}\label{equ101}\mathbf{v}:= \frac{\mathbf{u}}{\sqrt{\varepsilon(\Lambda)}}\sum_{i=0}^{m-1}\frac{\sqrt{\lambda(e_i)}}{\dim_{\mathbbm{k}}(V_i)}e_i.\end{equation} As we shall see, the element $\mathbf{v}$ plays a key role in the representation theory of smash product Hopf algebras. For the element $\mathbf{v}$, we have the following result:

\begin{prop}\label{prop2}The element $\mathbf{v}$ satisfies the the following properties:
\begin{enumerate}
  \item $\varepsilon(\mathbf{v})=1$.
  \item $S^2(h)=\mathbf{v}h\mathbf{v}^{-1}$ for $h\in H$.
  \item $\mathbf{v}^2=\mathbf{u}S(\mathbf{u}^{-1})$, which is the distinguished group-like element of $H$.
  \item $\mathbf{v}^n=1$, where $n=2\dim_{\mathbbm{k}}(H)$.
  \item $\mathbf{v}^{-1}=S(\mathbf{v})$.
  \item $\mathbf{v}=1$ if and only if $S^2=id$.
\end{enumerate}
\end{prop}
\proof (1) Applying $\varepsilon$ to both sides of the equality (\ref{equ101}), we obtain that $\varepsilon(\mathbf{v})=1$.

(2) Since $S^2(h)=\mathbf{u}h\mathbf{u}^{-1}$ and the elements $\mathbf{u}$, $\mathbf{v}$ are the same up to a central invertible element $\frac{1}{\sqrt{\varepsilon(\Lambda)}}\sum_{i=0}^{m-1}\frac{\sqrt{\lambda(e_i)}}{\dim_{\mathbbm{k}}(V_i)}e_i$, it follows that $S^2(h)=\mathbf{v}h\mathbf{v}^{-1}$ for $h\in H$.

(3) Note that $\mathbf{u}^{-1}S(\mathbf{u}^{-1})=\frac{1}{\varepsilon(\Lambda)}\sum_{i=0}^{m-1}\frac{\lambda(e_i)}{\dim_{\mathbbm{k}}(V_i)^2}e_i$ by Proposition \ref{prop3.20} (4). It follows that $$\mathbf{u}S(\mathbf{u}^{-1})=\frac{\mathbf{u}^2}{\varepsilon(\Lambda)}\sum_{i=0}^{m-1}\frac{\lambda(e_i)}{\dim_{\mathbbm{k}}(V_i)^2}e_i=\mathbf{v}^2,$$
which is the distinguished group-like element of $H$ by Proposition \ref{prop3.20} (5).

(4) It can be seen from Part (3) that $\mathbf{v}^2$ is the distinguished group-like element of $H$, while the order of the distinguished group-like element divides  $\dim_{\mathbbm{k}}(H)$. This implies that $\mathbf{v}^n=(\mathbf{v}^2)^{\dim_{\mathbbm{k}}(H)}=1$.

(5) There is a permutation $*$ on the index set $\{0,1,\cdots,m-1\}$ determined by $i^*=j$ if the dual $H$-module $V^*_{i}$ is isomorphic to $V_j$. The permutation $*$ satisfies that $i^{**}=i$, $S(e_i)=e_{i^*}$, $\dim_{\mathbbm{k}}(V_{i^*})=\dim_{\mathbbm{k}}(V_i)$ and $\lambda(e_{i^*})=\lambda(e_i)$ for $0\leq i\leq m-1$ (the last equality follows from \cite[Corollary 3.4]{WLL1}). We have
\begin{align*}
\mathbf{v}S(\mathbf{v})&=\frac{1}{\varepsilon(\Lambda)}\mathbf{u}\bigg(\sum_{i=0}^{m-1}\frac{\sqrt{\lambda(e_i)}}{\dim_{\mathbbm{k}}(V_i)}e_i\bigg)S(\mathbf{u})\bigg(\sum_{i=0}^{m-1}\frac{\sqrt{\lambda(e_i)}}{\dim_{\mathbbm{k}}(V_i)}e_{i^*}\bigg)\\
&=\frac{1}{\varepsilon(\Lambda)}\mathbf{u}S(\mathbf{u})\bigg(\sum_{i=0}^{m-1}\frac{\sqrt{\lambda(e_i)}}{\dim_{\mathbbm{k}}(V_i)}e_i\bigg)
\bigg(\sum_{i=0}^{m-1}\frac{\sqrt{\lambda(e_{i^*})}}{\dim_{\mathbbm{k}}(V_{i^*})}e_i\bigg)\\
&=\frac{1}{\varepsilon(\Lambda)}\mathbf{u}S(\mathbf{u})\bigg(\sum_{i=0}^{m-1}\frac{\sqrt{\lambda(e_i)}}{\dim_{\mathbbm{k}}(V_i)}e_i\bigg)^2\\
&=\mathbf{u}S(\mathbf{u})\bigg(\frac{1}{\varepsilon(\Lambda)}\sum_{i=0}^{m-1}\frac{\lambda(e_i)}{\dim_{\mathbbm{k}}(V_i)^2}e_i\bigg)\\
&=1,
\end{align*}
where the last equality follows from  Proposition \ref{prop3.20} (4). We obtain that $\mathbf{v}^{-1}=S(\mathbf{v})$.

(6) If $\mathbf{v}=1$, it follows from Part (2) that  $S^2=id$. Conversely, if $S^2=id$, then $\Lambda$ is cocommutative (see \cite[Proposition 3.5]{WLL1}) and hence $\mathbf{u}=\varepsilon(\Lambda)$. In this case, $\lambda(e_i)=\frac{\dim_{\mathbbm{k}}(V_i)^2}{\varepsilon(\Lambda)}$ by Proposition \ref{prop3.20} (3). Taking it into the equality (\ref{equ101}) we may see that $\mathbf{v}=1$.
We complete the proof.
\qed

\section{\bf Representations of smash product Hopf algebras}
We denote $n:=2\dim_{\mathbbm{k}}(H)$. Let $G$ be a cyclic group of order $n$ generated by $g$. The character group $\widehat{G}$ of $G$ is also a cyclic group of order $n$. Let $\psi$ be a generator of $\widehat{G}$. Then $\widehat{G}=\{\psi^j\mid 0\leq j\leq n-1\}$, which is the complete set of distinct irreducible characters of simple $\mathbbm{k}G$-modules. The simple $\mathbbm{k}G$-module with respect to the character $\psi^j$ is denoted by $W_j$ for $0\leq j\leq n-1$.

Since the antipode $S$ of $H$ satisfies $S^{2n}=id$ by Radford's formula of $S^4$ \cite{Rad1}, the Hopf algebra $H$ is a left $\mathbbm{k}G$-module algebra whose action is given by $$g^i\rightarrow h=S^{2i}(h)\ \text{for}\ g^i\in G\  \text{and}\ h\in H.$$
This reduces to a Hopf algebra $H\#\mathbbm{k}  G$ mentioned in \cite{So}. More precisely, the Hopf algebra $H\#\mathbbm{k}G$ is the smash product of $H$ and $\mathbbm{k}G$. The multiplication of $H\#\mathbbm{k}G$ is given by $$(a\# g^i)(b\#g^j)=a(g^i\rightarrow b)\# g^{i+j}=aS^{2i}(b)\# g^{i+j}\ \text{for}\ a,b\in H,$$ the identity of $H\#\mathbbm{k}G$ is $1_H\#1_{\mathbbm{k}  G}$. The comultiplication of $H\#\mathbbm{k}G$ is given by $$\Delta_{H\#\mathbbm{k}G}(h\# g^i)=(h_{(1)}\# g^i)\otimes(h_{(2)}\# g^i).$$ The counit of $H\#\mathbbm{k}G$ is $\varepsilon_{H\#\mathbbm{k}G}=\varepsilon_H\#\varepsilon_{\mathbbm{k}G}$ and the antipode of $H\#\mathbbm{k}G$ is $$S_{H\#\mathbbm{k}G}(h\# g^i)=(1_H\#g^{-i})(S(h)\#1_{\mathbbm{k}G})=S^{1-2i}(h)\# g^{-i}.$$ Moreover, $1_H\# g$ is a group-like element of $H\#\mathbbm{k} G$ that satisfies \begin{equation}\label{equ7-25}S_{H\#\mathbbm{k}G}^2(h\# g^i)=(1_H\# g)(h\# g^i)(1_H\# g)^{-1}.\end{equation} The Hopf algebra $H$ can be considered as a sub-Hopf algebra of $H\#\mathbbm{k}G$ under the injective map $H\rightarrow H\#\mathbbm{k}G,\ h\mapsto h\#1_{\mathbbm{k}G}$.

Since $\Lambda$ is an integral of $H$ with $\varepsilon(\Lambda)\neq0$ and $p\nmid n$, $\Lambda\#\frac{1}{n}\sum_{i=0}^{n-1}g^i$ is an integral of $H\#\mathbbm{k}G$ with $\varepsilon_{H\#\mathbbm{k}G}(\Lambda\#\frac{1}{n}\sum_{i=0}^{n-1}g^i)=\varepsilon(\Lambda)\neq0$. Thus, $H\#\mathbbm{k}G$ is a semisimple Hopf algebra over $\mathbbm{k}$.

The representation theory of crossed product of an algebra with a group algebra has been studied in \cite{MW}. However, we do not take advantage of those notations and methods in \cite{MW} to describe $H\#\mathbbm{k}G$-modules. Instead, since the Hopf algebra $H\#\mathbbm{k}G$ is semisimple, we will determine all simple $H\#\mathbbm{k}G$-modules by the study of the character of regular representation of $H\#\mathbbm{k}G$.

\begin{lem}
If $V$ is a finite dimensional $H$-module and $W$ is a finite dimensional $\mathbbm{k}G$-module, then the vector space $V\otimes W$ is a finite dimensional $H\#\mathbbm{k}G$-module, where the $H\#\mathbbm{k}G$-module structure is given by
\begin{equation}\label{equ000}(h\# g^k)\cdot(v\otimes w)=(h\mathbf{v}^k\cdot v)\otimes (g^k\cdot w)\ \text{for}\ v\in V, w\in W.\end{equation}
\end{lem}
\proof By Proposition \ref{prop2} (4), we have $\mathbf{v}^n=1$. It follows that
$$(h\# g^n)\cdot(v\otimes w)=(h\mathbf{v}^n\cdot v)\otimes (g^n\cdot w)=(h\cdot v)\otimes w=(h\# 1_{\mathbbm{k}G})\cdot(v\otimes w).$$ This is compatible with the equality $h\# g^n=h\# 1_{\mathbbm{k}G}$. For $a,b\in H$,
by $S^2(h)=\mathbf{v}h\mathbf{v}^{-1}$ for $h\in H$, we may check that
$$((a\# g^k)(b\# g^j))\cdot(v\otimes w)=(a\# g^k)\cdot((b\# g^j)\cdot(v\otimes w)).$$ The proof is completed.
\qed

\begin{lem}\label{lem2}
If $V$ is a simple $H$-module and $W$ is a simple  $\mathbbm{k}G$-module, then $V\otimes W$ is a simple $H\#\mathbbm{k}G$-module.
\end{lem}
\proof Note that $H\#\mathbbm{k}G$ is a semisimple Hopf algebra over an algebraically closed field $\mathbbm{k}$. It is sufficient to show that $\text{End}_{H\#\mathbbm{k}G}(V\otimes W)=\mathbbm{k}$. Suppose that the map $\phi:V\otimes W\rightarrow V\otimes W$ is an  $H\#\mathbbm{k}G$-module morphism. Since $W$ is one dimensional, we fix a basis $w$ of $W$. The $H\#\mathbbm{k}G$-module morphism $\phi$ induces an $H$-module morphism $\phi_0:V\rightarrow V$ as follows: $\phi(v\otimes w)=\phi_0(v)\otimes w$ for any $v\in V$.
This shows that $\phi$ is the identity map of $V\otimes W$ up to a scalar, since $V$ is simple and $\phi_0$ is the identity map of $V$ up to a scalar.
\qed

\begin{rem}For simple $H$-module $V_i$ and simple $\mathbbm{k}G$-module $W_j$, it can be seen from Lemma \ref{lem2} that $V_i\otimes W_j$ is a simple $H\#\mathbbm{k}G$-module. Let $\chi_{ij}$ be the character associated to the simple $H\#\mathbbm{k}G$-module $V_i\otimes W_j$. It follows from (\ref{equ000}) that $$\chi_{ij}(h\otimes g^k)=\chi_i(h\mathbf{v}^k)\psi^j(g^k)\ \text{for}\ 0\leq i\leq m-1, 0\leq j\leq n-1.$$
\end{rem}

\begin{thm}The set
$\{V_i\otimes W_j\mid 0\leq i\leq m-1,0\leq j\leq n-1\}$ forms a complete set of non-isomorphic simple $H\#\mathbbm{k}G$-modules.
\end{thm}
\proof
Note that $\Lambda\#\frac{1}{n}\sum_{i=0}^{n-1}g^i$ is a left integral in $H\#\mathbbm{k}G$ and $\lambda\#\sum_{j=0}^{n-1}\psi^j$ is a right integral in $(H\#\mathbbm{k}G)^*$ satisfying $(\lambda\#\sum_{j=0}^{n-1}\psi^j)(\Lambda\#\frac{1}{n}\sum_{i=0}^{n-1}g^i)=1.$ By \cite[Corollary 6]{Rad}, the characters of left regular representations of $H$ and $H\#\mathbbm{k}G$ are respectively given by $\chi_H=\lambda\leftharpoonup \mathbf{u}$ and $\chi_{H\#\mathbbm{k}G}=(\lambda\#\sum_{j=0}^{n-1}\psi^j)\leftharpoonup\mathbf{u}_{H\#\mathbbm{k}G}$, where $\mathbf{u}=S(\Lambda_{(2)})\Lambda_{(1)}$ and
\begin{align*}\mathbf{u}_{H\#\mathbbm{k}G}&=\frac{1}{n}\sum_{i=0}^{n-1} S_{H\#\mathbbm{k}  G}(\Lambda_{(2)}\#g^i)(\Lambda_{(1)}\#g^i)\\
&=\frac{1}{n}\sum_{i=0}^{n-1} (S^{1-2i}(\Lambda_{(2)})\# g^{-i})(\Lambda_{(1)}\# g^i)\\
&=\frac{1}{n}\sum_{i=0}^{n-1} S^{1-2i}(\Lambda_{(2)})S^{-2i}(\Lambda_{(1)})\# 1_{\mathbbm{k}G}\\
&=\frac{1}{n}\sum_{i=0}^{n-1} S^{-2i}(\mathbf{u})\# 1_{\mathbbm{k}G}\\
&=\mathbf{u}\# 1_{\mathbbm{k}G}.
\end{align*}
It follows that $$\chi_{H\#\mathbbm{k}G}=(\lambda\#\sum_{j=0}^{n-1}\psi^j)\leftharpoonup (\mathbf{u}\# 1_{\mathbbm{k}G})=(\lambda\leftharpoonup \mathbf{u})\#\sum_{j=0}^{n-1}\psi^j=\chi_H\#\sum_{j=0}^{n-1}\psi^j.$$
Hence, $$(\chi_{H\#\mathbbm{k}G})(h\#g^k)=\chi_H(h)\sum_{j=0}^{n-1}\psi^j(g^k)=\left\{
                                                                                 \begin{array}{ll}
                                                                                   n\chi_H(h), & k=0; \\
                                                                                   0, & 1\leq k\leq n-1.
                                                                                 \end{array}
                                                                               \right.
$$
While
\begin{align*}\sum_{i=0}^{m-1}\sum_{j=0}^{n-1}\dim_{\mathbbm{k}}(V_i\otimes W_j)\chi_{ij}(h\#g^k)&=\sum_{i=0}^{m-1}\sum_{j=0}^{n-1}\dim_{\mathbbm{k}}(V_i)\chi_i(h\mathbf{v}^k)\psi^j(g^k)\\
&=\chi_H(h\mathbf{v}^k)\sum_{j=0}^{n-1}\psi^j(g^k)\\
&=\left\{
                                                                                 \begin{array}{ll}
                                                                                   n\chi_H(h), & k=0; \\
                                                                                   0, & 1\leq k\leq n-1.
                                                                                 \end{array}
                                                                               \right.
\end{align*}
We obtain that
$\chi_{H\#\mathbbm{k}G}=\sum_{i=0}^{m-1}\sum_{j=0}^{n-1}\dim_{\mathbbm{k}}(V_i\otimes W_j)\chi_{ij}$. Hence, all non-isomorphic simple $H\#\mathbbm{k}G$-modules are $V_i\otimes W_j$ for $0\leq i\leq m-1,0\leq j\leq n-1$.
\qed

\begin{rem}Note that $\chi_{00}=\varepsilon_{H\#\mathbbm{k}G}.$ Hence $V_0\otimes W_0$ is
the trivial $H\#\mathbbm{k}G$-module, where $V_0$ is the trivial $H$-module and $W_0$ is the trivial $\mathbbm{k}G$-module.
\end{rem}

The dual module $(V_i\otimes W_j)^*$ of $V_i\otimes W_j$ can be described as follows:

\begin{prop}\label{p2}We have
$(V_i\otimes W_j)^*\cong V_{i^*}\otimes W_{j^*}$ for $0\leq i\leq m-1,0\leq j\leq n-1$, where $V_{i^*}$ is the dual of $V_i$ as an $H$-module and $W_{j^*}$ is the dual of $W_j$ as a $\mathbbm{k}G$-module.
\end{prop}
\proof We need to check that $\chi_{i^*j^*}=\chi_{ij}\circ S_{H\#\mathbbm{k}G}$ for $0\leq i\leq m-1,0\leq j\leq n-1$. Note that $S(\mathbf{v})=\mathbf{v}^{-1}$ and $S^{-2}(h)=\mathbf{v}^{-1}h\mathbf{v}$ for $h\in H$.
On the one hand,
$$\chi_{i^*j^*}(h\#g^k)=\chi_{i^*}(h\mathbf{v}^k)\psi^{j^*}(g^k)=\chi_{i}(S(\mathbf{v})^kS(h))\psi^j(g^{-k})=\chi_{i}(\mathbf{v}^{-k}S(h))\psi^j(g^{-k}).$$
On the other hand,
\begin{align*}(\chi_{ij}\circ S_{H\#\mathbbm{k}G})(h\#g^k)
&=\chi_{ij}(S_{H\#\mathbbm{k}G}(h\#g^k))\\
&=\chi_{ij}(S^{1-2k}(h)\#g^{-k})\\
&=\chi_{ij}(\mathbf{v}^{-k}S(h)\mathbf{v}^k\#g^{-k})\\
&=\chi_{i}(\mathbf{v}^{-k}S(h))\psi^{j}(g^{-k}).
\end{align*} We conclude that $\chi_{i^*j^*}=\chi_{ij}\circ S_{H\#\mathbbm{k}G}$ for $0\leq i\leq m-1,0\leq j\leq n-1$.
\qed

\section{\bf The Grothendieck algebras of smash product Hopf algebras}
In this section, we will investigate a relationship between the Grothendieck algebra of the smash product Hopf algebra $H\#\mathbbm{k}G$ and the Grothendieck algebra of $H$.
Recall that the Grothendieck algebra $(G_0(H)\otimes_{\mathbb{Z}}\mathbbm{k},*)$ of $H$ is an associative algebra over $\mathbbm{k}$ with unity $\varepsilon_H$ under the convolution $*$, where the convolution $*$ is defined by
$$(\chi_i\ast \chi_j)(h)=(\chi_i\otimes \chi_j)(\Delta(h))\ \text{for}\ h\in H.$$
%The set $\{\chi_i\mid 0\leq i\leq m-1\}$ forms a $\mathbbm{k}$-basis of $G_0(H)$.
We define a new multiplication operator $\star$ on $G_0(H)\otimes_{\mathbb{Z}}\mathbbm{k}$ by
$$(\chi_i\star \chi_j)(h)=(\chi_i\otimes \chi_j)\bigg(\Delta(h)\Delta(\mathbf{v}^{-1})(\mathbf{v}\otimes \mathbf{v})\bigg)\ \text{for}\ h\in H.$$

\begin{rem}\label{rem2}
If $S^2=id,$ then $\mathbf{v}=1$ by Proposition \ref{prop2} (6). In this case, the multiplication operator $\star$ is nothing but the convolution $\ast$.
\end{rem}

\begin{prop}
The pair $(G_0(H)\otimes_{\mathbb{Z}}\mathbbm{k},\star)$ is an associative algebra over $\mathbbm{k}$ with unity $\varepsilon_H$.
\end{prop}
\proof
We first prove that $\star$ is a multiplication operator on $G_0(H)\otimes_{\mathbb{Z}}\mathbbm{k}$. That is, $\chi_i\star\chi_j\in G_0(H)\otimes_{\mathbb{Z}}\mathbbm{k}$ for $0\leq i,j\leq m-1$.
Indeed, for $a,b\in H$, using $S^2(h)=\mathbf{v}h\mathbf{v}^{-1}$ for $h\in H$, we have
\begin{align*}(\chi_i\star\chi_j)(ab)&=\chi_i(a_{(1)}b_{(1)}{\mathbf{v}^{-1}}_{(1)}\mathbf{v})\chi_j(a_{(2)}b_{(2)}{\mathbf{v}^{-1}}_{(2)}\mathbf{v})\\
&=\chi_i(a_{(1)}(b\mathbf{v}^{-1})_{(1)}\mathbf{v})\chi_j(a_{(2)}(b\mathbf{v}^{-1})_{(2)}\mathbf{v})\\
&=\chi_i(a_{(1)}(\mathbf{v}^{-1}S^2(b))_{(1)}\mathbf{v})\chi_j(a_{(2)}(\mathbf{v}^{-1}S^2(b))_{(2)}\mathbf{v})\\
&=\chi_i(a_{(1)}{\mathbf{v}^{-1}}_{(1)}S^2(b_{(1)})\mathbf{v})\chi_j(a_{(2)}{\mathbf{v}^{-1}}_{(2)}S^2(b_{(2)})\mathbf{v})\\
&=\chi_i(a_{(1)}{\mathbf{v}^{-1}}_{(1)}\mathbf{v}b_{(1)})\chi_j(a_{(2)}{\mathbf{v}^{-1}}_{(2)}\mathbf{v}b_{(2)})\\
&=\chi_i(b_{(1)}a_{(1)}{\mathbf{v}^{-1}}_{(1)}\mathbf{v})\chi_j(b_{(2)}a_{(2)}{\mathbf{v}^{-1}}_{(2)}\mathbf{v})\\
&=(\chi_i\star\chi_j)(ba).
\end{align*}
It follows from \cite{Lo} that  $\chi_i\star\chi_j\in G_0(H)\otimes_{\mathbb{Z}}\mathbbm{k}$ for $0\leq i,j\leq m-1$. Since the map $H\rightarrow H\otimes H,\ h\mapsto\Delta(h)\Delta(\mathbf{v}^{-1})(\mathbf{v}\otimes \mathbf{v})$ is a coassociative comultiplication in $H$ for which $\varepsilon_H$ is still a counit (see \cite[Eq.(12)]{AEGN}), the operator $\star$ dual to the coassociative comultiplication is an associative multiplication on $G_0(H)\otimes_{\mathbb{Z}}\mathbbm{k}$ with unity $\varepsilon_H$.
\qed
%The associativity and unity $\varepsilon$ of $\star$ on $G_0(H)\otimes_{\mathbb{Z}}\mathbbm{k}$ can be checked directly. Indeed, for $a\in H$, we have
%\begin{align*}((\chi_i\star\chi_j)\star\chi_k)(a)&=((\chi_i\star\chi_j)\otimes\chi_k)\bigg(\Delta(a)\Delta(\mathbf{v}^{-1})(\mathbf{v}\otimes \mathbf{v})\bigg)\\
%&=(\chi_i\star\chi_j)(a_{(1)}{\mathbf{v}^{-1}}_{(1)}\mathbf{v})\chi_k(a_{(2)}{\mathbf{v}^{-1}}_{(2)}\mathbf{v})\\
%&=\chi_i(a_{(1)}{\mathbf{v}^{-1}}_{(1)}\mathbf{v})\chi_j(a_{(2)}{\mathbf{v}^{-1}}_{(2)}\mathbf{v})\chi_k(a_{(3)}{\mathbf{v}^{-1}}_{(3)}\mathbf{v})\\
%&=\chi_i(a_{(1)}{\mathbf{v}^{-1}}_{(1)}\mathbf{v})(\chi_j\star\chi_k)(a_{(2)}{\mathbf{v}^{-1}}_{(2)}\mathbf{v})\\
%&=(\chi_i\star(\chi_j\star\chi_k))(a).
%\end{align*}
%Therefore, $((\chi_i\star\chi_j)\star\chi_k)=(\chi_i\star(\chi_j\star\chi_k))$ for $0\leq i,j,k\leq m-1$.
%\begin{align*}(\varepsilon\star\chi_k)(a)&=(\varepsilon\otimes\chi_k)\bigg(\Delta(a)\Delta(\mathbf{v}^{-1})(\mathbf{v}\otimes \mathbf{v})\bigg)\\
%&=\varepsilon(a_{(1)}{\mathbf{v}^{-1}}_{(1)}\mathbf{v})\chi_k(a_{(2)}{\mathbf{v}^{-1}}_{(2)}\mathbf{v})\\
%&=\chi_k(a).
%\end{align*}Hence, $\varepsilon\star\chi_k=\chi_k$ for $0\leq k\leq m-1$.It is similar that $\chi_k\star\varepsilon=\chi_k$ for $0\leq k\leq m-1$.

Next, we will use the algebras $(G_0(H)\otimes_{\mathbb{Z}}\mathbbm{k},*)$ and $(G_0(H)\otimes_{\mathbb{Z}}\mathbbm{k},\star)$ to describe the structure of the Grothendieck algebra $(G_0(H\#\mathbbm{k}G)\otimes_{\mathbb{Z}}\mathbbm{k},\ast)$ of $H\#\mathbbm{k}G$. Note that $\{\chi_0,\chi_1,\cdots,\chi_{m-1}\}$ is a $\mathbbm{k}$-basis of  $G_0(H)\otimes_{\mathbb{Z}}\mathbbm{k}$.
Suppose in $(G_0(H)\otimes_{\mathbb{Z}}\mathbbm{k},*)$ that
$$\chi_i\ast\chi_j=\sum_{k=0}^{n-1}N_{ij}^k\chi_k$$
and in $(G_0(H)\otimes_{\mathbb{Z}}\mathbbm{k},\star)$ that $$\chi_i\star\chi_j=\sum_{k=0}^{n-1}L_{ij}^k\chi_k,$$
where $N_{ij}^k$ and $L_{ij}^k$ are respectively the structure coefficients of the two algebras with respect to the basis $\{\chi_0,\chi_1,\cdots,\chi_{m-1}\}$. We stress that the coefficient $N_{ij}^k$ is the multiplicity of $V_k$ appeared in the decomposition of tensor product $V_i\otimes V_j$ as $H$-modules, so each $N_{ij}^k$  is indeed a nonnegative integer. For the coefficient $L_{ij}^k$, we shall see in Remark \ref{rem3} that each $L_{ij}^k$ is an integer.
%For the the Grothendieck algebra $G_0(H\#\mathbbm{k}G)\otimes_{\mathbb{Z}}\mathbbm{k}$ of $H\#\mathbbm{k}G$, we have the following result:
\begin{prop}\label{prop1}
We have the following equations in the Grothendieck algebra $(G_0(H\#\mathbbm{k}G)\otimes_{\mathbb{Z}}\mathbbm{k},\ast)$:
\begin{enumerate}
  \item $\chi_{ij}=\chi_{i0}\ast\chi_{0j}=\chi_{0j}\ast\chi_{i0}$ for $0\leq i\leq m-1,0\leq j\leq n-1$.
  \item $\chi_{i0}\ast\chi_{j0}=\sum_{k=0}^{m-1}\frac{1}{2}(N_{ij}^k+L_{ij}^k)\chi_{k0}+\sum_{k=0}^{m-1}\frac{1}{2}(N_{ij}^k-L_{ij}^k)\chi_{k\frac{n}{2}}$ for $0\leq i,j\leq m-1$.
  \item $\chi_{is}\ast\chi_{jt}=\sum_{k=0}^{m-1}\frac{1}{2}(N_{ij}^k+L_{ij}^k)\chi_{k,s+t}+\sum_{k=0}^{m-1}\frac{1}{2}(N_{ij}^k-L_{ij}^k)\chi_{k,\frac{n}{2}+s+t}$ for $0\leq i,j\leq m-1$ and $0\leq s,t\leq n-1$, where $s+t$ and $\frac{n}{2}+s+t$ are reduced modulo $n$.
\end{enumerate}
\end{prop}
\proof (1) It is direct to calculate that
\begin{align*}(\chi_{i0}\ast\chi_{0j})(h\#g^k)&=\chi_{i0}(h_{(1)}\#g^k)\chi_{0j}(h_{(2)}\#g^k)\\
&=\chi_i(h_{(1)}\mathbf{v}^k)\psi^0(g^k)\chi_0(h_{(2)}\mathbf{v}^k)\psi^j(g^k)\\
&=\chi_i(h\mathbf{v}^k)\psi^j(g^k)\\
&=\chi_{ij}(h\#g^k).
\end{align*}
So we have $\chi_{i0}\ast\chi_{0j}=\chi_{ij}$. It is similar that $\chi_{0j}\ast\chi_{i0}=\chi_{ij}$.

%Note that $\mathbf{v}^2$ is the distinguished group-like element of $H$. We have
%\begin{align*}
%(\chi_{i0}\ast\chi_{j0})(a\#g^{2s})&=\chi_{i0}(a_{(1)}\#g^{2s})\chi_{j0}(a_{(2)}\#g^{2s})\\
%&=\chi_{i}(a_{(1)}\mathbf{v}^{2s})\chi_{j}(a_{(2)}\mathbf{v}^{2s})\\
%&=(\chi_{i}*\chi_{j})(a\mathbf{v}^{2s})\\
%&=\sum_{k=0}^{m-1}N_{ij}^k\chi_k(a\mathbf{v}^{2s}).
%\end{align*}

(2) We show that the values that both sides of the desired equation taking on $h\#g^{l}$ are the same. Note that $\mathbf{v}^2$ is the distinguished group-like element of $H$ and $\psi^{\frac{n}{2}}(g)=-1$.
For the case $l=2s$, we have
\begin{align*}
&\ \ \ \ \sum_{k=0}^{m-1}\frac{1}{2}(N_{ij}^k+L_{ij}^k)\chi_{k0}(h\#g^{2s})+\sum_{k=0}^{m-1}\frac{1}{2}(N_{ij}^k-L_{ij}^k)\chi_{k\frac{n}{2}}(h\#g^{2s})\\
&=\sum_{k=0}^{m-1}\frac{1}{2}(N_{ij}^k+L_{ij}^k)\chi_{k}(h\mathbf{v}^{2s})+\sum_{k=0}^{m-1}\frac{1}{2}(N_{ij}^k-L_{ij}^k)\chi_{k}(h\mathbf{v}^{2s})\psi^{\frac{n}{2}}(g^{2s})\\
&=\sum_{k=0}^{m-1}\frac{1}{2}(N_{ij}^k+L_{ij}^k)\chi_{k}(h\mathbf{v}^{2s})+\sum_{k=0}^{m-1}\frac{1}{2}(N_{ij}^k-L_{ij}^k)\chi_{k}(h\mathbf{v}^{2s})\\
&=\sum_{k=0}^{m-1}N_{ij}^k\chi_{k}(h\mathbf{v}^{2s})=(\chi_{i}*\chi_{j})(h\mathbf{v}^{2s})\\
&=\chi_{i}(h_{(1)}\mathbf{v}^{2s})\chi_{j}(h_{(2)}\mathbf{v}^{2s})\ \ \ \ (\text{since}\ \mathbf{v}^{2s}\ \text{is\ a\ group-like\ element})\\
&=\chi_{i0}(h_{(1)}\#g^{2s})\chi_{j0}(h_{(2)}\#g^{2s})=(\chi_{i0}\ast\chi_{j0})(h\#g^{2s}).
\end{align*}
For the case $l=2s+1$, we have
\begin{align*}
&\ \ \ \ \sum_{k=0}^{m-1}\frac{1}{2}(N_{ij}^k+L_{ij}^k)\chi_{k0}(h\#g^{2s+1})+\sum_{k=0}^{m-1}\frac{1}{2}(N_{ij}^k-L_{ij}^k)\chi_{k\frac{n}{2}}(h\#g^{2s+1})\\
&=\sum_{k=0}^{m-1}\frac{1}{2}(N_{ij}^k+L_{ij}^k)\chi_{k}(h\mathbf{v}^{2s+1})+\sum_{k=0}^{m-1}\frac{1}{2}(N_{ij}^k-L_{ij}^k)\chi_{k}(h\mathbf{v}^{2s+1})\psi^{\frac{n}{2}}(g^{2s+1})\\
&=\sum_{k=0}^{m-1}\frac{1}{2}(N_{ij}^k+L_{ij}^k)\chi_{k}(h\mathbf{v}^{2s+1})-\sum_{k=0}^{m-1}\frac{1}{2}(N_{ij}^k-L_{ij}^k)\chi_{k}(h\mathbf{v}^{2s+1})\\
&=\sum_{k=0}^{m-1}L_{ij}^k\chi_{k}(h\mathbf{v}^{2s+1})=(\chi_{i}\star\chi_{j})(h\mathbf{v}^{2s+1})\\
&=\chi_{i}(h_{(1)}\mathbf{v}^{2s+1})\chi_{j}(h_{(2)}\mathbf{v}^{2s+1})\ \ \ \ (\text{since}\ \mathbf{v}^{2s}\ \text{is\ a\ group-like\ element})\\
&=\chi_{i0}(h_{(1)}\#g^{2s+1})\chi_{j0}(h_{(2)}\#g^{2s+1})=(\chi_{i0}\ast\chi_{j0})(h\#g^{2s+1}).
\end{align*}
We obtain the desired equation.

(3) Using Part (1) and Part (2) we may see that Part (3) holds.
\qed

\begin{rem}\label{rem3} It follows from Proposition \ref{prop1} (2) that the tensor product $(V_i\otimes W_0)\otimes(V_j\otimes W_0)$ has the following decomposition as $H\#\mathbbm{k}G$-modules:
$$(V_i\otimes W_0)\otimes(V_j\otimes W_0)\cong\bigoplus_{k=0}^{m-1}\frac{1}{2}(N_{ij}^k+L_{ij}^k)(V_{k}\otimes W_0)\bigoplus\bigoplus_{k=0}^{m-1}\frac{1}{2}(N_{ij}^k-L_{ij}^k)(V_k\otimes W_{\frac{n}{2}}).$$
Thus, these coefficients $\frac{1}{2}(N_{ij}^k+L_{ij}^k)$ and $\frac{1}{2}(N_{ij}^k-L_{ij}^k)$ are both nonnegative integers. Since all $N^k_{ij}$ are nonnegative integers, it follows that all $L^k_{ij}$ are integers and satisfy $-N^k_{ij}\leq L^k_{ij}\leq N^k_{ij}$. In view of this, the multiplication operator $\star$ defined on the Grothendieck algebra $G_0(H)\otimes_{\mathbb{Z}}\mathbbm{k}$ can be defined as well on the Grothendieck ring $G_0(H).$
\end{rem}

The Grothendieck algebra $(G_0(H\#\mathbbm{k}G)\otimes_{\mathbb{Z}}\mathbbm{k},\ast)$ is an associative unity algebra with a $\mathbbm{k}$-basis $\{\chi_{ij}\mid 0\leq i\leq m-1,0\leq j\leq n-1\}$.
Denote by $$\theta_l=\frac{1}{n}\sum_{t=0}^{n-1}\psi(g)^{-lt}\chi_{0t}\ \text{for}\ 0\leq l\leq n-1.$$ Note that $\chi_{0t}=\psi^t$ for $0\leq t\leq n-1$. Thus, $\{\theta_l\mid0\leq l\leq n-1\}$ is the set of all central primitive idempotents of the algebra $\mathbbm{k}\widehat{G}$.
Moreover, we have \begin{equation}\label{equ0}\chi_{0j}*\theta_{l}=\psi(g)^{jl}\theta_l\ \text{and}\ \chi_{ij}*\theta_{l}=\chi_{i0}*\chi_{0j}*\theta_{l}=\psi(g)^{jl}\chi_{i0}*\theta_l.\end{equation} In particular, each $\theta_l$ is a central idempotent of $(G_0(H\#\mathbbm{k}G)\otimes_{\mathbb{Z}}\mathbbm{k},\ast)$. The structure of the Grothendieck algebra $(G_0(H\#\mathbbm{k}G)\otimes_{\mathbb{Z}}\mathbbm{k},\ast)$ now can be described as follows:

%\begin{prop}For any $0\leq l\leq n-1$,  the element $\theta_l$ satisfies the following properties:
%\begin{enumerate}
%  \item $\chi_{0j}*\theta_{l}=\psi(g)^{jl}\theta_l$.
%  \item $\chi_{ij}*\theta_{l}=\chi_{i0}*\chi_{0j}*\theta_{l}=\psi(g)^{jl}\chi_{i0}*\theta_l$.
%  \item $\theta_l$ is a central idempotent of $G_0(H\#\mathbbm{k}G)\otimes_{\mathbb{Z}}\mathbbm{k}$.
%\end{enumerate}
%\end{prop}

\begin{thm}\label{th1}We have the following algebra isomorphisms:
\begin{enumerate}
  \item If $l$ is even, then $(G_0(H\#\mathbbm{k}G)\otimes_{\mathbb{Z}}\mathbbm{k},\ast)*\theta_l\cong(G_0(H)\otimes_{\mathbb{Z}}\mathbbm{k},*)$.
  \item If $l$ is odd, then $(G_0(H\#\mathbbm{k}G)\otimes_{\mathbb{Z}}\mathbbm{k},\ast)*\theta_l\cong(G_0(H)\otimes_{\mathbb{Z}}\mathbbm{k},\star)$.
  \item We have $(G_0(H\#\mathbbm{k}G)\otimes_{\mathbb{Z}}\mathbbm{k},\ast)\cong(G_0(H)\otimes_{\mathbb{Z}}\mathbbm{k},*)^{\oplus\frac{n}{2}}
\bigoplus(G_0(H)\otimes_{\mathbb{Z}}\mathbbm{k},\star)^{\oplus\frac{n}{2}}.$
\end{enumerate}
\end{thm}
\proof (1) For the case $l$ being even, we consider the $\mathbbm{k}$-linear map $$\phi_l:(G_0(H)\otimes_{\mathbb{Z}}\mathbbm{k},*)\rightarrow (G_0(H\#\mathbbm{k}G)\otimes_{\mathbb{Z}}\mathbbm{k},\ast)*\theta_l,\ \ \chi_i\mapsto\chi_{i0}*\theta_l.$$
It can be seen from (\ref{equ0}) that $\phi_l$ is bijective, and moreover, $\chi_{i\frac{n}{2}}\ast\theta_l=\chi_{i0}\ast\theta_l$. Now
\begin{align*}\phi_l(\chi_i*\chi_j)&=\phi_l(\sum_{k=0}^{m-1}N_{ij}^k\chi_k)=\sum_{k=0}^{m-1}N_{ij}^k\chi_{k0}*\theta_l\\
&=\sum_{k=0}^{m-1}\frac{1}{2}(N_{ij}^k+L_{ij}^k)\chi_{k0}*\theta_l+\sum_{k=0}^{m-1}\frac{1}{2}(N_{ij}^k-L_{ij}^k)\chi_{k0}*\theta_l\\
&=\sum_{k=0}^{m-1}\frac{1}{2}(N_{ij}^k+L_{ij}^k)\chi_{k0}*\theta_l+\sum_{k=0}^{m-1}\frac{1}{2}(N_{ij}^k-L_{ij}^k)\chi_{k\frac{n}{2}}*\theta_l\\
&=(\chi_{i0}*\chi_{j0})*\theta_l=(\chi_{i0}*\theta_l)\ast(\chi_{j0}*\theta_l)\\
&=\phi_l(\chi_i)*\phi_l(\chi_j).
\end{align*}
This shows that $\phi_l$ is an algebra isomorphism.

(2) For the case $l$ being odd, we consider the $\mathbbm{k}$-linear map $$\phi_l:(G_0(H)\otimes_{\mathbb{Z}}\mathbbm{k},\star)\rightarrow (G_0(H\#\mathbbm{k}G)\otimes_{\mathbb{Z}}\mathbbm{k},\ast)*\theta_l,\ \ \chi_i\mapsto\chi_{i0}*\theta_l.$$
It can be seen from (\ref{equ0}) that $\phi_l$ is bijective, and moreover, $\chi_{i\frac{n}{2}}\ast\theta_l=-\chi_{i0}\ast\theta_l$. Now
\begin{align*}\phi_l(\chi_i\star\chi_j)&=\phi_l(\sum_{k=0}^{m-1}L_{ij}^k\chi_k)=\sum_{k=0}^{m-1}L_{ij}^k\chi_{k0}*\theta_l\\
&=\sum_{k=0}^{m-1}\frac{1}{2}(N_{ij}^k+L_{ij}^k)\chi_{k0}*\theta_l-\sum_{k=0}^{m-1}\frac{1}{2}(N_{ij}^k-L_{ij}^k)\chi_{k0}*\theta_l\\
&=\sum_{k=0}^{m-1}\frac{1}{2}(N_{ij}^k+L_{ij}^k)\chi_{k0}*\theta_l+\sum_{k=0}^{m-1}\frac{1}{2}(N_{ij}^k-L_{ij}^k)\chi_{k\frac{n}{2}}*\theta_l\\
&=(\chi_{i0}*\chi_{j0})*\theta_l=(\chi_{i0}*\theta_l)\ast(\chi_{j0}*\theta_l)\\
&=\phi_l(\chi_i)*\phi_l(\chi_j).
\end{align*}
Thus, $\phi_l$ is an algebra isomorphism.

(3) Let $(G_0(H)\otimes_{\mathbb{Z}}\mathbbm{k},*)^{\oplus\frac{n}{2}}$ be the direct sum of $\frac{n}{2}$-folds of $(G_0(H)\otimes_{\mathbb{Z}}\mathbbm{k},*)$ and $(G_0(H)\otimes_{\mathbb{Z}}\mathbbm{k},\star)^{\oplus\frac{n}{2}}$ the direct sum of $\frac{n}{2}$-folds of $(G_0(H)\otimes_{\mathbb{Z}}\mathbbm{k},\star)$. Since $\theta_0+\theta_1+\cdots+\theta_{n-1}=1$, where $1$ is the identity $\chi_{00}$ of  $(G_0(H\#\mathbbm{k}G)\otimes_{\mathbb{Z}}\mathbbm{k},\ast)$, using Part (1) and Part (2) we obtain the following algebra isomorphism: $$(G_0(H\#\mathbbm{k}G)\otimes_{\mathbb{Z}}\mathbbm{k},\ast)\cong(G_0(H)\otimes_{\mathbb{Z}}\mathbbm{k},*)^{\oplus\frac{n}{2}}
\bigoplus(G_0(H)\otimes_{\mathbb{Z}}\mathbbm{k},\star)^{\oplus\frac{n}{2}}.$$
The proof is completed. \qed

\begin{rem}
If $S^2=id$, by Remark \ref{rem2}, the algebra $(G_0(H)\otimes_{\mathbb{Z}}\mathbbm{k},\star)$ is nothing but the Grothendieck algebra $(G_0(H)\otimes_{\mathbb{Z}}\mathbbm{k},*)$. In this case,
$$(G_0(H\#\mathbbm{k}G)\otimes_{\mathbb{Z}}\mathbbm{k},\ast)\cong(G_0(H)\otimes_{\mathbb{Z}}\mathbbm{k},*)^{\oplus n}.$$
\end{rem}

Let $\mathcal{C}$ be the $\mathbbm{k}$-linear subcategory of Rep$(H\#\mathbbm{k}G)$ spanned by objects $$\{V_i\otimes W_0,V_i\otimes W_{\frac{n}{2}}\mid 0\leq i\leq m-1\}.$$ Then $\mathcal{C}$ is closed under taking dual by Proposition \ref{p2}. It follows from Proposition \ref{prop1} that $\mathcal{C}$ is also closed under the tensor product of objects. More explicitly,
$$(V_i\otimes W_{\frac{n}{2}})\otimes(V_j\otimes W_{\frac{n}{2}})\cong(V_i\otimes W_0)\otimes(V_j\otimes W_0)$$$$\cong\bigoplus_{k=0}^{m-1}\frac{1}{2}(N_{ij}^k+L_{ij}^k)(V_{k}\otimes W_0)\bigoplus\bigoplus_{k=0}^{m-1}\frac{1}{2}(N_{ij}^k-L_{ij}^k)(V_k\otimes W_{\frac{n}{2}}),$$
and
$$(V_i\otimes W_0)\otimes(V_j\otimes W_{\frac{n}{2}})\cong(V_i\otimes W_{\frac{n}{2}})\otimes(V_j\otimes W_{0})$$$$\cong\bigoplus_{k=0}^{m-1}\frac{1}{2}(N_{ij}^k+L_{ij}^k)(V_{k}\otimes W_{\frac{n}{2}})\bigoplus\bigoplus_{k=0}^{m-1}\frac{1}{2}(N_{ij}^k-L_{ij}^k)(V_k\otimes W_{0}).$$
Hence $\mathcal{C}$ is a fusion subcategory of Rep$(H\#\mathbbm{k}G)$.
Let $(G_0(\mathcal{C})\otimes_{\mathbb{Z}}\mathbbm{k},\ast)$ be the Grothendieck algebra of $\mathcal{C}$. Then $\{\chi_{i0},\chi_{i\frac{n}{2}}\mid0\leq i\leq m-1\}$ forms a $\mathbbm{k}$-basis of $(G_0(\mathcal{C})\otimes_{\mathbb{Z}}\mathbbm{k},\ast)$.

\begin{prop}\label{p1}
We have the following algebra isomorphism: $$(G_0(\mathcal{C})\otimes_{\mathbb{Z}}\mathbbm{k},\ast)\cong (G_0(H)\otimes_{\mathbb{Z}}\mathbbm{k},*)\bigoplus(G_0(H)\otimes_{\mathbb{Z}}\mathbbm{k},\star).$$
\end{prop}
\proof We denote $\theta=\frac{1}{2}(\chi_{00}+\chi_{0\frac{n}{2}})$. Then $1-\theta=\frac{1}{2}(\chi_{00}-\chi_{0\frac{n}{2}}),$ where $1$ is the identity $\chi_{00}$ of  $(G_0(\mathcal{C})\otimes_{\mathbb{Z}}\mathbbm{k},\ast)$. Note that
$\theta$ and $1-\theta$ are both central idempotents of $(G_0(\mathcal{C})\otimes_{\mathbb{Z}}\mathbbm{k},\ast)$. In particular, $$\chi_{i\frac{n}{2}}*\theta=\chi_{i0}*\chi_{0\frac{n}{2}}*\theta=\chi_{i0}*\theta\ \text{for}\ 0\leq i\leq m-1.$$
Consider the $\mathbbm{k}$-linear map $$\phi:(G_0(H)\otimes_{\mathbb{Z}}\mathbbm{k},*)\rightarrow (G_0(\mathcal{C})\otimes_{\mathbb{Z}}\mathbbm{k},\ast)*\theta,\ \ \chi_i\mapsto\chi_{i0}*\theta.$$
It is easy to see that $\phi$ is bijective and
\begin{align*}\phi(\chi_i*\chi_j)&=\phi(\sum_{k=0}^{m-1}N_{ij}^k\chi_k)=\sum_{k=0}^{m-1}N_{ij}^k\chi_{k0}*\theta\\
&=\sum_{k=0}^{m-1}\frac{1}{2}(N_{ij}^k+L_{ij}^k)\chi_{k0}*\theta+\sum_{k=0}^{m-1}\frac{1}{2}(N_{ij}^k-L_{ij}^k)\chi_{k0}*\theta\\
&=\sum_{k=0}^{m-1}\frac{1}{2}(N_{ij}^k+L_{ij}^k)\chi_{k0}*\theta+\sum_{k=0}^{m-1}\frac{1}{2}(N_{ij}^k-L_{ij}^k)\chi_{k\frac{n}{2}}*\theta\\
&=(\chi_{i0}*\chi_{j0})*\theta=(\chi_{i0}*\theta)\ast(\chi_{j0}*\theta)\\
&=\phi(\chi_i)*\phi(\chi_j).
\end{align*}
This shows that $\phi$ is an algebra isomorphism.
Consider the $\mathbbm{k}$-linear map $$\varphi:(G_0(H)\otimes_{\mathbb{Z}}\mathbbm{k},\star)\rightarrow (G_0(\mathcal{C})\otimes_{\mathbb{Z}}\mathbbm{k},\ast)*(1-\theta),\ \ \chi_i\mapsto\chi_{i0}*(1-\theta).$$
Then $\varphi$ is bijective. Using $\chi_{i\frac{n}{2}}*(1-\theta)=-\chi_{i0}*(1-\theta)$ we may see that
\begin{align*}\varphi(\chi_i\star\chi_j)&=\varphi(\sum_{k=0}^{m-1}L_{ij}^k\chi_k)=\sum_{k=0}^{m-1}L_{ij}^k\chi_{k0}*(1-\theta)\\
&=\sum_{k=0}^{m-1}\frac{1}{2}(N_{ij}^k+L_{ij}^k)\chi_{k0}*(1-\theta)-\sum_{k=0}^{m-1}\frac{1}{2}(N_{ij}^k-L_{ij}^k)\chi_{k0}*(1-\theta)\\
&=\sum_{k=0}^{m-1}\frac{1}{2}(N_{ij}^k+L_{ij}^k)\chi_{k0}*(1-\theta)+\sum_{k=0}^{m-1}\frac{1}{2}(N_{ij}^k-L_{ij}^k)\chi_{k\frac{n}{2}}*(1-\theta)\\
&=(\chi_{i0}*\chi_{j0})*(1-\theta)=(\chi_{i0}*(1-\theta))\ast(\chi_{j0}*(1-\theta))\\
&=\varphi(\chi_i)*\varphi(\chi_j)
\end{align*}
Hence,
$\varphi$ is an algebra isomorphism.
\qed

Note that $\theta=\theta_0+\theta_2+\theta_4+\cdots+\theta_{n-2}$ and $1-\theta=\theta_1+\theta_3+\theta_5+\cdots+\theta_{n-1}$. By Theorem \ref{th1} and Proposition \ref{p1}, we the following corollary:
\begin{cor}
We have algebra isomorphism: $$(G_0(H\#\mathbbm{k}G)\otimes_{\mathbb{Z}}\mathbbm{k},\ast)\cong(G_0(\mathcal{C})\otimes_{\mathbb{Z}}\mathbbm{k},\ast)^{\oplus\frac{n}{2}}.$$
\end{cor}

%The higher Frobenius-Schur indicators of finite dimensional representations of semisimple Hopf algebras over a field of characteristic 0 have been studied in \cite{KSZ}. In this section, we will consider these indicators of finite dimensional representations of the semisimple Hopf algebra $H\#\mathbbm{k}G$  over the field $\mathbbm{k}$ of characteristic $p>\dim_{\mathbbm{k}}(H)^{1/2}$.

Finally, we give some remarks on the pivotal (spherical) structure of the fusion categories Rep$(H\#\mathbbm{k}G)$ and $\mathcal{C}$.
Since $S^2_{H\#\mathbbm{k}G}$ is an inner automorphism of $H\#\mathbbm{k}G$ and $$S^2_{H\#\mathbbm{k}G}(h\# g^i)=(1_H\# g)(h\# g^i)(1_H\# g)^{-1},$$ where $1_H\# g$ is a group-like element of $H\#\mathbbm{k}G$, the category Rep$(H\#\mathbbm{k}G)$ is a pivotal fusion category, where the pivotal structure $\tau$ on Rep$(H\#\mathbbm{k}G)$ is the isomorphism of monoidal functors $\tau_{V\otimes W}:V\otimes W\rightarrow (V\otimes W)^{**}$ natural in $V\otimes W$. Here $\tau_{V\otimes W}(v\otimes w)$ is defined by $$\tau_{V\otimes W}(v\otimes w)(f)=f(1_H\# g\cdot v\otimes w)=f(\mathbf{v}\cdot v\otimes g\cdot w)$$ for $v\in V, w\in W$ and $f\in(V\otimes W)^{*}$.

The quantum dimension of $V\otimes W\in \text{Rep}(H\#\mathbbm{k}G)$ with respect to the pivotal structure $\tau$ is denoted by $\mathbf{dim}(V\otimes W)$, which is the following composition
$$\mathbf{1}\xrightarrow{\text{coev}_{(V\otimes W)}}(V\otimes W)\otimes(V\otimes W)^*\xrightarrow{\tau_{V\otimes W}\otimes id}(V\otimes W)^{**}\otimes(V\otimes W)^*\xrightarrow{\text{ev}_{(V\otimes W)^*}}\mathbf{1},$$where $\mathbf{1}$ is the trivial $H\#\mathbbm{k}G$-module $V_0\otimes W_0$. From this composition, we have
$$\mathbf{dim}(V\otimes W)=\chi_V(\mathbf{v})\chi_W(g).$$
Especially,
$$\mathbf{dim}(V_i\otimes W_j)=\chi_i(\mathbf{v})\psi^j(g)=\sqrt{\varepsilon(\Lambda)}\sqrt{\lambda(e_i)}\psi^j(g).$$ For the dual module $(V_i\otimes W_j)^*\cong V_{i^*}\otimes W_{j^*}$, we have $$\mathbf{dim}(V_{i^*}\otimes W_{j^*})=\sqrt{\varepsilon(\Lambda)}\sqrt{\lambda(e_{i^*})}\psi^j(g^{-1})=\sqrt{\varepsilon(\Lambda)}\sqrt{\lambda(e_{i})}\psi^j(g^{-1}).$$
Therefore, $\mathbf{dim}(V_{i^*}\otimes W_{j^*})=\mathbf{dim}(V_{i}\otimes W_{j})$ if and only if $j=0$ or $j=\frac{n}{2}$. This means that  with respect to the pivotal structure $\tau$, the fusion category
Rep$(H\#\mathbbm{k}G)$ is pivotal but not spherical, while the the fusion subcategory $\mathcal{C}$ of Rep$(H\#\mathbbm{k}G)$ spanned by objects $\{V_i\otimes W_0,V_i\otimes W_{\frac{n}{2}}\mid 0\leq i\leq m-1\}$ is both pivotal and spherical.

{\bf Question:} Does the pivotalization of Rep$(H)$ equal to the pivotal fusion category $\mathcal{C}$?

\section*{Acknowledgement}
The first author was supported by Qing Lan Project.
The second author was supported by National Natural Science Foundation of China (Grant No. 12171230). The third author was supported by National Natural Science Foundation of China (Grant No. 11871063).


\vskip5pt
\begin{thebibliography}{99}
\bibliographystyle{siam}

%\bibitem{An}
%N. Andruskiewitsch, Notes on extensions of Hopf algebras, Can. J. Math., 48 (1), 1996: 3-42.


%\bibitem{BFS}
%K. I. Beidar, Y. Fong, A. A. Stolin, On antipodes and integrals in Hopf algebras over rings and the quantum Yang-Baxter equation, \textit{J. Algebra} {\bf194}(1) (1997) 36-52.

%\bibitem{Cib}
%C. Cibils, A quiver quantum group, Communications in mathematical physics, 1993, 157(3): 459-477.

%\bibitem{CW}
%M. Cohen, S. Westreich, Probabilistically nilpotent Hopf algebras, Transactions of the American Mathematical Society, 368(6), 2016, 4295-4314.

%\bibitem{Doi}
%Y. Doi, Group-like algebras and their representations, Communications in Algebra, {\bf 38} (7), (2010): 2635-2655.

%\bibitem{ENO}
%P. Etingof, D. Nikshych, V. Ostrik, On fusion categories, Annals of Mathematics, 2005, 162, 581-642.


%\bibitem{GN}
%S. Gelaki, D. Nikshych, Nilpotent fusion categories, Adv. Math. 217, 1053-1071 (2008).


%\bibitem{LR1}R. G. Larson, D. E. Radford, Semisimple cosemisimple Hopf algebras, \textit{Am. J. Math.}, {\bf109} (1987) 187-195.

%\bibitem{LR}R. G. Larson, D. E. Radford, Finite dimensional cosemisimple Hopf algebras in characteristic 0 are semisimple, \textit{J. Algebra} {\bf117} (1988) 267-289.
\bibitem{AEGN}
E. Aljadeff, P. Etingof, S. Gelaki, D. Nikshych, On twisting of finite-dimensional Hopf algebras, J. Algebra {\bf256} (2002): 484-501.

\bibitem{CYW}
J. Chen, S. Yang, D. Wang, Grothendieck rings of a class of Hopf algebras of Kac-Paljutkin type, Front. Math. China {\bf16} (1) (2021): 29-47.

\bibitem{EG}
P. Etingof, S. Gelaki, On finite dimensional semisimple and cosemisimple Hopf algebras in positive characteristic, Internat. Math. Res. Notices (16) (1998): 851-864.


\bibitem{Ka}
Y. Kashina, Classification of semisimple Hopf algebras of dimension 16, J. Algebra {\bf232} (2000): 617-663.

\bibitem{LL}
K. Li, G. Liu, On the antipode of Hopf algebras with the dual Chevalley property, J. Pure Appl. Algebra {\bf226} (3) (2022): 106871.

\bibitem{Lo}
M. Lorenz, Representations of finite-dimensional Hopf algebras, J. Algebra {\bf 188} (1997): 476-505.

\bibitem{Mo}
R. K. Molnar, Semi-direct products of Hopf algebras, J. Algebra {\bf47} (1) (1977): 29-51.

\bibitem{Mon}
S. Montgomery, Hopf Algebras and their actions on rings, CBMS Series in Math., Vol. {\bf 82}, AMS, Providence, 1993.

\bibitem{MW}
S. Montgomery, S. J. Witherspoon, Irreducible representations of crossed products, J. Pure Appl. Algebra {\bf129} (3) (1998): 315-326.

\bibitem{NR}
W. D. Nichols, M. B. Richmond, The Grothendieck group of a Hopf algebra, J. Pure Appl. Algebra {\bf106}  (3) (1996): 297-306.

\bibitem{Ni}
D. Nikshych, $K_0$-rings and twisting of finite dimensional semisimple Hopf algebras, Comm. Algebra {\bf26} (1998): 321-342.

\bibitem{Rad}
D. E. Radford, The trace function and Hopf algebras, J. Algebra {\bf163} (1994): 583-622.

\bibitem{Rad1}
D. E. Radford, The order of the antipode of a finite dimensional Hopf algebra is finite, Amer. J. Math. {\bf 98} (1976): 333-355.


\bibitem{So}
Y. Sommerh$\ddot{a}$user, On Kaplansky's fifth conjecture, J. Algebra {\bf 204} (1998): 202-224.


\bibitem{WLL1}
Z. Wang, G. Liu, L. Li, Higher Frobenius-Schur indicators for semisimple Hopf algebras in positive characteristic, arXiv:2112.15264v2.


\bibitem{YY}
R. Yang, S. Yang, The Grothendieck rings of Wu-Liu-Ding algebras and their Casimir numbers (II), Comm.  Algebra {\bf49} (5) (2021): 2041-2073.



%\bibitem{Ng}
%S.-H. Ng, A Note on Frobenius-Schur Indicators, Proc. of the International Conference on Algebra 2010, World Sci. Publ., Hackensack, NJ (2012),454-460.


%\bibitem{GP}
%M. Geck, G. Pfeiffer, Characters of finite Coxeter groups and Iwahori-Hecke algebras, New York: Oxford University Press, 2000.

%\bibitem{Lar}
%R. G. Larson, Characters of Hopf algebras, J. Algebra 17 (1971), 352-368.

%\bibitem{Lor}
%M. Lorenz, Some applications of Frobenius algebras to Hopf algebras, Contemp. Math., 2011,  537, 269-289.


%\bibitem{Ch}
%Ryba Christopher. The Structure of the Grothendieck Rings of Wreath Product Deligne Categories and their Generalisations[J]. International Mathematics Research Notices,2019,2021(16):

%\bibitem{}
%Seung-moon Hong,Eric Rowell. On the classification of the Grothendieck rings of non-self-dual modular categories[J]. Journal of Algebra,2010,324(5):
%\bibitem{Mon}
%S. Montgomery, Hopf Algebras and their actions on rings, CBMS Series in Math., Vol. {\bf 82}, AMS, Providence, 1993.


\end{thebibliography}
\end{document}